\documentclass{article}%
\usepackage{amsmath}%
\usepackage{amsfonts}%
\usepackage{amssymb}%
\usepackage{graphicx}%
\usepackage{color} 

\newtheorem{thm}{Theorem}

\newtheorem{cor}[thm]{Corollary}
\newtheorem{pro}[thm]{Proposition}

\def\a{\alpha }
\def\b{\beta }

\def\t{\theta}
\def\g{\gamma }

\def\Bbb{\mathbb}

\def\tr{\mathop{\rm tr}}

\def\re{\mathop{\rm {Re}}}

\def\wm<{\prec_{\tiny w}}
\def\wm>{\succ_{\tiny w}}
\def\lm<{\prec_{\mbox{\tiny $\log$}}}
\def\lm>{\succ_{\mbox{\tiny $\log$}}}
\def\wl<{\prec_{\mbox{\rm \tiny w$\log$}}}
\def\wl>{\succ_{\mbox{\rm \tiny w$\log$}}}

\def\psd<{\prec}
\def\psd>{\succ}
\def\psdgeq>{\succeq}

\def\n<{<}
\def\n>{>}

\def\proof{{\noindent {\bf Proof.} \hspace{.01in}}}
\newcommand{\qed}{\Box \smallskip}

\begin{document}


\title{\Large Angles, triangle inequalities, correlation matrices and  metric-preserving and subadditive functions\thanks{to appear in Linear Algebra Appl.}}

\author{Diego Castano\thanks{E-mails: castanod@nova.edu},
Vehbi E. Paksoy\thanks{E-mails: vp80@nova.edu}, Fuzhen Zhang\thanks{E-mails: zhang@nova.edu}\\
\\
\small{Farquhar College of Arts and Sciences,
Nova Southeastern University}\\
\small{3301 College Ave, Fort Lauderdale, FL 33314, USA}}
\date{}
\maketitle

\hrule
\bigskip
\noindent {\bf Abstract}

\medskip
We present inequalities concerning the entries of correlation matrices, density matrices, and partial isometries
through the positivity of $3\times 3$ matrices.  We extend our discussions to the inequalities
concerning the triangle triplets with metric-preserving and subadditive functions.

\medskip
{\footnotesize

}

\bigskip
\hrule





\section{Introduction}

We begin our discussions with angles between vectors by looking at
 two angle definitions and  their characteristics. We
 first study the triangle inequalities for the angles through the positivity (i.e., positive semidefiniteness)
of a $3\times 3$ matrix; we investigate the relations between the triangle inequalities
(of angles or more generally the triangle triplets) and metric-preserving and subadditive
  functions.
Our results will capture some existing ones but through a different approach.
The discussion on the $3\times 3$ positive semidefinite matrices leads to some inequalities
concerning the entries of correlation matrices with which we obtain inequalities for density matrices and
partial isometries.

Let $V$ be an inner product space with the inner product $\langle \cdot  , \cdot \rangle $
over the real number field $\Bbb R$.
For  any two nonzero vectors $u, v$ in $V$, there are two common ways to
define the {\it angle} between the vectors $u$ and $v$ in terms of the inner product (see, for instance,
 \cite[p.\,58]{Lau05} and \cite[p.\,335]{Lay},  respectively):
\begin{eqnarray}
\theta \left( u, v\right) &=&\arccos \frac{\left\vert \left\langle
u,v\right\rangle \right\vert}{\|u\|\, \|v\|}  
\label{eqt}
\\
\Theta \left( u,v \right) &=&\arccos \frac{\left\langle u, v\right\rangle}{\|u\|\, \|v\|}
\label{eqT}
\end{eqnarray}%

There are various reasons that the angles are defined in ways (\ref{eqt}) and (\ref{eqT})
(in the sense of Euclidean geometry).
Definition (\ref{eqt}) may stem from the angles between subspaces, while (\ref{eqT})
makes perfect sense intuitively. We are interested in the properties of the angles regardless their definitions.

Angle and inner product can be viewed as a ``twin" for a vector space.
 A vector may take a simple and familiar form like the ones in $\Bbb R^2$; it may look much more complicated
like the elements (linear combinations of wedge products) in the Grassmann spaces \cite[p.172]{MerMul97}. Some matrix functions
 are closely related to vectors and some types of  products of vectors.
It is a well-known fact that the trace of a matrix product is an inner product:
$\tr (AB)=\langle B, A^*\rangle$. The determinant (even more generally,
the generalized matrix functions) can be expressed as an inner product of $*$-tensors
(see, e.g., \cite[p.\,226]{MerMul97}).

We will also need the term {\em correlation matrix}, which is a positive semidefinite matrix with all main diagonal
entries equal to 1. Every positive semidefinite matrix with nonzero main  diagonal entries  can be normalized to a
 correlation matrix through scaling. The correlation  matrices are frequently used in statistics.
 For its determinant and permanent
properties, see, e.g.,  \cite{Pie87, Wat88, FZ13Per}. Our theorems rely on the results
for  the  $3\times 3$ correlation matrices.

  In Section 2, we focus on the triangle inequalities through the positivity of $3\times 3$ matrices. Our results provide
  a unified proof for the triangle inequalities for the angles $\t$ and $\Theta$. We also present some relationship
   between the elements of correlation matrices. As applications, we obtain inequalities for
   density matrices and partial isometries. In Section 3, we study the triangle triplets
  (which are more general than the angles formed by three vectors),  metric-preserving and subadditive functions. Some inequalities of unit
  vectors are immediate from our results.

\section{The triangle inequality and $3\times 3$ matrices}

We extend somewhat the underlying number field of the vector space to the complex number field $\Bbb C$
(see \cite[p.\,9]{DK70})   and replace the angle in (\ref{eqT}) by
\begin{eqnarray}
\Theta \left( u,v \right) &=&\arccos \frac{\re \left\langle u, v\right\rangle}{\|u\|\, \|v\|}
\label{eqT2}
\end{eqnarray}%

The angles $\t$ and $\Theta$ are closely related, but not equal unless $\left\langle u, v\right\rangle$ is nonnegative.
Since $\re \left\langle u, v\right\rangle\leq \left\vert \left\langle
u,v\right\rangle \right\vert$ and $f(t)=\arccos t$ is a decreasing function in $t\in [-1, 1]$, we have $\Theta \geq \theta$.
On the other hand, if 
$\left\langle v, u\right\rangle\neq 0$,
by taking $p = \frac{\left\langle v, u\right\rangle}{|\left\langle v, u\right\rangle|}$, we get
$\theta (u, v)=\Theta (p u, v).$  It is easy to verify that  (see, e.g.,  \cite{Lin12_MathIntel})
\begin{eqnarray}
\theta \left( u, v\right) = \min_{|p|=1}\Theta \left(p u,v \right) =
\min_{|q|=1}\Theta \left( u,q v \right)=\min_{|p|=|q |=1}\Theta \left(pu, qv  \right)
\label{eqt2}
\end{eqnarray}%

For any nonzero vectors $u, v\in V$, 
we see $\theta \left( u, v\right)\in [0, \frac{\pi}{2}]$,
$\Theta \left( u, v\right)\in [0, \pi]$, and $\theta (u, -u)=0$,
while $\Theta (u, -u)=\pi$. The angle $\theta$ defined in (\ref{eqt})
between $u$ and $v$ is $\frac{\pi}{2}$ if and only if $u$ and $v$
are {\em orthogonal}, i.e., $\langle u, v\rangle =0$; however,
the angles $\theta$ in general  do not obey the law of cosines
(for the triangle formed by
nonzero vectors $u, v$ and $u-v$).
 In contrast, the law of cosines does hold for the
 angles $\Theta$, 
 and $\Theta \left( u, v\right)=\frac{\pi}{2}$ if and only if
 $\langle u, v\rangle =0$  over $\Bbb R$,    but it is possible for some  vectors $u, v$
 to have
 an angle $\Theta (u, v) =\frac{\pi}{2}$ and $\langle u, v\rangle \neq 0$  over $\Bbb C$.

The triangle inequalities for $\theta $ and $\Theta$ are known. That is,
for all nonzero vectors $u, v, w\in V$ and the angles $\theta $ in (\ref{eqt}) and  $\Theta$
in (\ref{eqT2}),

\begin{itemize}
\setlength{\itemindent}{.2in}
\item[(t).]\label{zt}  $\theta (u, v) \leq \theta (u, w)+\theta (w, v)$.
\item[(T).]\label{zT}  $\Theta (u, v) \leq \Theta (u, w)+\Theta (w, v)$.
\end{itemize}

The triangle inequality (T) is attributed to Krein \cite{Kre69} by Gustafson and Rao \cite[p.\,56]{GR97}.
The inequality was stated without proof in \cite{Kre69} and proved first in \cite{DKR76}, then  in
\cite[p.\,56]{GR97}. Note that the real case for (T) is also seen in \cite[p.\,31]{ZFZbook11}.
It has been observed \cite{Lin12_MathIntel} that (t) follows from (T) because of
(\ref{eqt2}).

The proof of (T) in \cite[p.\,56]{GR97} boils down to the positivity of the matrix
$$R_0=\left ( \begin{array}{ccc}
  1 & \re \langle u, v\rangle & \re\langle u, w\rangle \\
  \re\langle v, u\rangle & 1 & \re\langle v, w\rangle \\
 \re\langle w, u\rangle & \re\langle w, v\rangle & 1
 \end{array} \right )
 $$\mbox{for unit vectors $u, v, w$,} which is ensured by the positivity of the Gram matrix
$$G_0=\left ( \begin{array}{ccc}
  1 & \langle u, v\rangle & \langle u, w\rangle \\
  \langle v, u\rangle & 1 & \langle v, w\rangle \\
 \langle w, u\rangle & \langle w, v\rangle & 1
 \end{array} \right )$$
The positivity of $G_0$ also guarantees (see, e.g., \cite[p.\,26]{BhaPDM07}) the positivity  of
 $$
A_0= \left ( \begin{array}{ccc}
  1 & |\langle u, v\rangle| & |\langle u, w\rangle| \\
 | \langle v, u\rangle| & 1 & |\langle v, w\rangle| \\
 |\langle w, u\rangle| & |\langle w, v\rangle| & 1
 \end{array} \right )$$
which results in
\begin{equation}\label{Eq1}
1+2|\langle u, v\rangle |\,|\langle v, w\rangle |\,|\langle w, u\rangle |\geq
          |\langle u, v\rangle |^2+|\langle v, w\rangle |^2+|\langle w, u\rangle |^2
          \end{equation}
Inequality (\ref{Eq1}) is weaker than the following existing  inequality (see (5.2) of Theorem 5.1 in \cite{ZFZActa12}):
\begin{equation}\label{Eq1b}
1+2\re \big (\langle u, v\rangle \langle v, w\rangle \langle w, u\rangle \big )\geq
          |\langle u, v\rangle |^2+|\langle v, w\rangle |^2+|\langle w, u\rangle |^2
  \end{equation}

Using the idea of Gustafson and Rao
 we present a unified proof for the triangle inequalities (t) and (T) through $3\times 3$
matrices. We have known that
$3\times 3$ matrices play important roles in geometry and analysis. For instance,
the area of a triangle and the volume of a parallelepiped
in $\Bbb R^3$, as well as
the convexity of real-valued functions,  can be computed and determined by
(or through the determinants of) $3\times 3$ matrices (see, e.g., \cite[p.\,2]{PPT92}).

\begin{pro}\label{Thm1:arccos}
Let $a, b, c$ be real numbers such that the $3\times 3$ matrix
$$B=\left ( \begin{array}{ccc}
  1 & a & b  \\
  a  & 1 & c \\
b & c   & 1
 \end{array} \right )$$
is positive semidefinite. Let $f(t)$ be a strictly  decreasing function defined on an interval
$\Bbb I\subseteq [0, \infty)$
with range  $[-1, 1]$.  
If
 \begin{equation}\label{arc}
 f(p+q)=f(p)f(q)-\sqrt{1-f^2(p)}\, \sqrt{1-f^2(q)}
 \end{equation}
 for all $p$, $q$, and $p+q$ in  $\Bbb I$,
 then for $x, y, z$ of  any arrangement of $a, b, c$,
 $$f^{-1}(x)\leq f^{-1}(y)+f^{-1}(z)$$
\end{pro}

\proof Since $B$ is positive semidefinite, we have $-1\leq x, y, z \leq 1$ and
$$1+2xyz\geq x^2+y^2+z^2$$
in which $x, y, z$ is any arrangement of $a, b, c$.
The above inequality implies
$$(1-y^2)(1-z^2)\geq x^2-2xyz+(yz)^2=(x-yz)^2$$
and
$$\sqrt{1-y^2}\sqrt{1-z^2}\geq |x-yz|$$
Hence
\begin{equation}\label{xyz}
x\geq yz-\sqrt{1-y^2}\sqrt{1-z^2}
\end{equation}

Note that as $f$ is decreasing, so is $f^{-1}$.
If $f^{-1}(-1)\leq f^{-1}(y)+f^{-1}(z)$, then $f^{-1}(x)\leq f^{-1}(-1)\leq f^{-1}(y)+f^{-1}(z)$.
Otherwise, $f^{-1}(-1)>f^{-1}(y)+f^{-1}(z)\geq f^{-1}(1)+0=f^{-1}(1)$; so,
$f^{-1}(y)+f^{-1}(z)$ lies in the domain of $f$.
Applying (\ref{arc}), we rewrite (\ref{xyz})   as
$$f(f^{-1}(x))\geq f(f^{-1}(y)+f^{-1}(z))$$ 
Since $f(x)$ is decreasing, we have $f^{-1}(x)\leq f^{-1}(y)+f^{-1}(z)$, as desired.
$\qed$

One may verify that
every function $f_{r}(t)=\cos (rt)$ with $r>0$ on $[0, \pi/r]$  satisfies
the condition (\ref{arc}).

\begin{cor}\label{cor3.14} The following statements follow from Proposition \ref{Thm1:arccos} immediately.
\begin{itemize}
\item[\rm (i)]  For any nonzero vectors $u,$ $v,$ and $w$, inequalities (t) and (T) hold.
\item[\rm (ii)]  Let $0\leq a, b, c\leq 1$ and  denote $\alpha =\arccos a$, $\b=\arccos b,$ $\g=\arccos c$. Then
all of the three inequalities $\a \leq \b+\g$, $\b\leq \g + \a$,
 and $ \g \leq \a +\b$   hold if and only if the matrix $B$ in Proposition \ref{Thm1:arccos} is positive
semidefinite.
\end{itemize}
\end{cor}

\proof (i). In Proposition \ref{Thm1:arccos},  setting $f(t)=\cos t$ on $[0, \pi]$ and taking the $3\times 3$ matrix $B$
to be the matrices $A_0$ and $R_0$, we obtain
the triangle inequalities (t) and (T), respectively. For (ii), one direction is clear. For the other direction,
observe that $\a\leq \b +\g\leq \pi$ when $0\leq a, b, c\leq 1$ and that
$0\leq \b-\g \leq \a$ or $0\leq \g-\b \leq \a$.
Because cosine  is decreasing on $[0, \pi]$,
we have
$bc-\sqrt{1-b^2}\sqrt{1-c^2}\leq a\leq  bc+\sqrt{1-b^2}\sqrt{1-c^2}$. It follows that
$1+2abc \geq a^2+b^2+c^2$. Hence, matrix $B$ is positive semidefinite.   
$\qed$

Note that for $(a, b, c)=(-1, -1, -1)$, we have
$(\arccos a, \arccos b, \arccos c)=(\pi, \pi, \pi)$,  the matrix $B$
in Proposition \ref{Thm1:arccos}
is not positive semidefinite.
Therefore, it is necessary to assume that $a, b, c$  are nonnegative
in Corollary \ref{cor3.14} (ii). Note that Corollary \ref{cor3.14} (ii)
has appeared in a recent paper \cite[Proposition 1.4]{Drury14} with a different proof.

We now present some inequalities concerning the entries of correlation matrices. Although the following theorem is stated
for the matrices of size $n\times n$, one can see that in essence it is  really about $3\times 3$ matrices.

\begin{thm}\label{main}
Let $A=(a_{ij})$ be an $n\times n$ complex correlation matrix.
Then for all integers $1\leq i, j, p, q\leq n$ such that $i<p<q$, $i< j<q$, and   any real $k\geq 2$,
$$\big |\,|a_{ip}|-|a_{iq}|\, \big |\leq \sqrt{1-|a_{jq}|^2}\leq \sqrt{2}\cdot  \sqrt{1-|a_{jq}|}$$
$$\big |\, |a_{ip}|^k-|a_{iq}|^k\big |\leq \sqrt{1-|a_{jq}|^k}$$

Analogous results hold  for $\re (a_{st})$ in place of  $a_{st}$ for all $s, t$.
\end{thm}
\proof This is immediate from   Theorem \ref{Thm:piday} and   Corollary \ref{cor1ref} below. 
$\qed$

Note that every principal submatrix of a positive semidefinite matrix is again positive semidefinite and   that
if $M=(m_{ij})$ is a $3\times 3$ positive semidefinite matrix, then so is
$N=(|m_{ij}|)$. (This is not true for $4\times 4$ or higher dimensions.)

Density matrices play an important role in quantum  computation; density matrices are the positive semidefinite matrices having trace 1.
\begin{cor}\label{quantum}
Let $G_i$ be density matrices and let $G_i=H_i^*H_i$,
where $H_i$ are $n\times n$ matrices, $i=1, 2, \dots, n$. Then for all $i<p<q$, $i< j<q$, and real $k\geq 2$,
$$\big |\,|\tr H^*_{i}H_p|-|\tr H^*_{i}H_q|\, \big |\leq \sqrt{1-|\tr H^*_{j}H_q|^2}
\leq \sqrt{2}\cdot  \sqrt{1-|\tr H^*_{j}H_q |}$$
$$\big |\, |\tr H^*_{i}H_p|^k-|\tr H^*_{i}H_q|^k\big |\leq \sqrt{1-|\tr H^*_{j}H_q|^k}$$

Analogous results hold  for $\re (\tr H^*_{s}H_t)$ and also for $\tr |H^*_{s}H_t|$ in place of  $\tr H^*_{s}H_t $ for all $s, t$.
Here, for a matrix $X$, $|X|$ stands for $(X^*X)^{1/2}$.
\end{cor}

\proof If the partitioned matrix $A=(A_{ij})$ is positive semidefinite,  then so is   $(\tr A_{ij})$ (see, e.g., \cite{ZFZActa12}).
Note that the matrix $(\tr H^*_iH_j)$ is correlation. The assertion for $ \tr |H^*_{s}H_t|$  is due to a recent result of Drury \cite{Drury14} (see also \cite{LiZ14}).
$\qed$

\begin{cor}\label{quantum}
Let $H_1, H_2, \dots, H_n$ be partial isometries each having $n$ columns,
i.e., every $H_i^*H_i=I_n$.
Then for all $i<p<q$, $i< j<q$, and real $k\geq 2$,
$$\big |\,|\det H^*_{i}H_p|-|\det H^*_{i}H_q|\, \big |\leq \sqrt{1-|\det H^*_{j}H_q|^2}\leq \sqrt{2}\cdot  \sqrt{1-|\det H^*_{j}H_q|}$$
$$\big | \, |\det H^*_{i}H_p|^k-|\det H^*_{i}H_q|^k\big |\leq \sqrt{1-|\det H^*_{j}H_q|^k}$$

Analogous results hold for
$\re (\det H^*_{s}H_t)$ in place of  $\det H^*_{s}H_t $ for all $s, t$.
\end{cor}

\proof If the partitioned matrix $A=(A_{ij})$ is positive semidefinite,  then so is  $(\det A_{ij})$ (see, e.g., \cite{ZFZActa12}).
Note that the matrix $(\det H^*_iH_j)$ is correlation. $\qed$

The proof of Theorem \ref{main}
reduces
to the following results on  $3\times 3$ correlation matrices. These results can be stated in terms of complex matrices for which
$a, b, c$ in the inequalities are replaced by $\re a, \re b, \re c$ or $|a|, |b|, |c|$, respectively.

\begin{thm}\label{Thm:piday}
Let $a, b, c$ be real numbers such that the $3\times 3$ matrix
$$B=\left ( \begin{array}{ccc}
  1 & a & b  \\
  a  & 1 & c \\
b & c   & 1
 \end{array} \right )$$
 is positive semidefinite.
 Denote
 $$c_-=ab-\sqrt{(1-a^2)(1-b^2)}, \quad c_+=ab+\sqrt{(1-a^2)(1-b^2)}$$
 and  let
 $$\Delta_{a, b}=
 \max  \Big \{\sqrt{1-c_-^2}, \sqrt{1-c_+^2} \,\Big \}, \;\;\delta_{a, b}=\min \Big \{\sqrt{1-c_-^2}, \sqrt{1-c_+^2}\, \Big \}$$
\begin{itemize}
\item[\rm (i).] If $f(x)$ is a function defined on the interval  $[c_{-},\, c_{+}]$ such that
 $$\delta_{a,b} \leq  f(x), \; \mbox{for all $x\in [c_{-}, \, c_{+}]$}$$ then
\begin{equation}\label{squarez}
|a^2-b^2| \leq \Delta_{a,b} f(x), \; \mbox{for all $x\in [c_{-},\,  c_{+}]$}
\end{equation}
 In particular
 \begin{equation}\label{square}
 |a^2-b^2|\leq \Delta_{a,b} \sqrt{1-c^{2}}
 \end{equation}
 \item[\rm (ii).]
 If $g(x)$ is a function defined on the interval  $[c_{-},\, c_{+}]$ such that
 $$\sqrt{1-c_+}\leq  g(x), \; \mbox{for all $x\in [c_{-}, \, c_{+}]$}$$ then
 \begin{equation}\label{one}
 |a-b|\leq \sqrt{1-c_-}\cdot  g(x), \; \mbox{for all $x\in [c_{-}, \, c_{+}]$}
 \end{equation}
  In particular
  \begin{equation}\label{one2}
 |a-b|\leq \sqrt{1-c_-} \cdot \sqrt{1-c}
 \end{equation}
\item[\rm (iii).]  If $a, b, c$ are in $[0, 1]$, then
 \begin{equation}\label{one}
 |a-b|\leq \sqrt{1-c^2}
 \end{equation}
 \end{itemize}
 \end{thm}

 \proof
 (i). Observe that $B$ is positive semidefinite if and only if $a, b, c\in [-1, 1]$ and
 $1+2abc\geq a^2+b^2+c^2$ and that $1+2abc\geq a^2+b^2+c^2$ if and only if
 $m(c):=c^2-2ab c -(1-a^2-b^2)\leq 0$. Note that $m(c)=0$ has
  solutions $c_-$, $c_+.$
We see that
  $m(x)\leq 0$ 
   for all $x\in [c_{-},\,  c_{+}]$. Moreover,  since $B$ is positive semidefinite,
   the scalar $c$ in the matrix $B$ lies in $[c_{-},\,  c_{+}]$.

 Let $a=\cos \mu$, $b=\cos \nu$,
 $\mu, \nu \in [0, \pi]$. Then $c_{\pm}=\cos (\mu\mp \nu)$. We compute
\begin{eqnarray*}
|a^2-b^2| & = & |\cos ^2 \mu -\cos ^2 \nu |
   =    |\sin (\mu+\nu)\sin (\mu-\nu)| \\
   & = &  \sqrt{1-c_-^2} \cdot  \sqrt{1-c_+^2} \\
  & = &  \max \Big \{\sqrt{1-c_-^2}, \sqrt{1-c_+^2} \,\Big  \}\cdot \min \Big \{\sqrt{1-c_-^2}, \sqrt{1-c_+^2} \, \Big  \}\\
  & \leq & \Delta_{a,b} f(x)
 \end{eqnarray*}
which is inequality (\ref{squarez}). Taking $f(x)=\sqrt{1-x^2}$, we have $f(x)\geq \delta_{a,b}$ for all $x \in [c_{-},\,  c_{+}]$.
 This   leads to  (\ref{square}) by setting $x=c$.

(ii). In a similar way by using trigonometric identities, we have
 \begin{eqnarray*}
|a-b| & = & |\cos  \mu -\cos \nu |
   =    \sqrt{1-\cos (\mu+\nu)} \cdot \sqrt{1-\cos(\mu-\nu)}\\
   & = &  \sqrt{1-c_-} \cdot  \sqrt{1-c_+} \\
  & \leq & \sqrt{1-c_-}  \cdot g(x)
 \end{eqnarray*}
 The special case
 (\ref{one2}) is because  $\sqrt{1-c_+}\leq \sqrt{1-c}$ for $c\in [c_-, c_+]$.

 (iii). If $0\leq a, b, c\leq 1$, then  $1+2abc \geq a^2+b^2+c^2$ yields
 $$|a-b|^2\leq -c^2+2abc+1-2ab=1-c^2+2ab(c-1)\leq 1-c^2$$  
Thus inequality (\ref{one}) follows. $\qed$

Note that $y=1+\frac{1}{c_{+}}\left (\sqrt{1-c_+^2}-1\right )x$ with $c_{+}\neq 0$
 is also a
function satisfying the conditions (i) and (ii) in Theorem \ref{Thm:piday}.
Equalities in (\ref{square}) and (\ref{one2}) occur when $B$ is the positive
semidefinite matrix with $a=1$, $b=c=0$. So in this sense the upper bounds for these inequalities are optimal.

We point out that   the restriction on $a, b, c$ in (\ref{one}) being nonnegative cannot be removed. 
For instance,
take $a=1$, 
$b=-1$ 
and $c=-1$. Then  $B$  is positive
semidefinite. However, $|a-b|> \sqrt{1-c^2}$.
Moreover, inequality (\ref{square})   reveals
$$|a^2-b^2|\leq   \sqrt{1-c^{2}}$$
as a sister inequality of (\ref{one}). For $k\geq 3$, $|a^k-b^k|$ is not bounded by $\sqrt{1-c^2}$ in general.
One may verify by the following example   that $|a^3-b^3|> \sqrt{1-c^2}$.
 {\rm  Let
$$C=\left ( \begin{array}{rrr}
  1  & 1    & 0.1  \\
  1  & 1    & 0.1 \\
0.1  & 0.1  & 1
 \end{array} \right )$$
 Then $C$ is positive semidefinite. For $a=1, b=c=0.1$, we have
 $$|a-b|=0.9, \;\; |a^2-b^2|=0.99, \;\; |a^3-b^3|=0.999, \;\; \sqrt{1-c^2}\approx 0.995$$}

For $x, y, z$ of any arrangement of  $a, b, c$ in the matrix $B$, Theorem \ref{Thm:piday} gives
\begin{equation}\label{Eq:13xyz}
|x^2-y^2|\leq \sqrt{1-z^2}
\end{equation}
and
\begin{equation}\label{Eq:13xyz2}
| x-y| \leq \sqrt{2}\cdot \sqrt{1-z}
\end{equation}

Inequalities (\ref{one}) and (\ref{Eq:13xyz}) 
imply that 
for any unit vectors $u, v, w$ and $k=1, 2$,
 \begin{equation*}
\Big | | \langle u, v\rangle |^k - |\langle u, w\rangle |^k \Big | \leq \sqrt{1-|\langle w, v\rangle |^2}
 \end{equation*}

Note that the positivity of $B$ in the previous theorem is
equivalent to $1+2abc\geq a^2+b^2+c^2$ for real numbers $a, b, c$ in $[-1, 1]$,
%
in which  $a, b, c$  are symmetric.

\begin{cor}\label{cor1ref}
Let $a, b, c$ be real numbers such that the $3\times 3$ matrix
$$B=\left ( \begin{array}{ccc}
  1 & a & b  \\
  a  & 1 & c \\
b & c   & 1
 \end{array} \right )$$
 is positive semidefinite.  Then the following statements hold.
 \begin{itemize}
 \item[\rm (i).]$\big | |a|-|b| \big |\leq  \sqrt{1-|c|^2}\leq \sqrt{2}\cdot \sqrt{1-|c|}$
     and for real $k\geq 2$,
 $$
 \big | |a|^k-|b|^k \big |\leq \sqrt{1-|c|^{k}}$$
  \item[\rm (ii).] If $0\leq \a, \b, \g \leq \pi/2$ and $|\a-\b |\leq \g \leq \a+\b$, then for any integer  $k\geq 1$
  $$|\cos ^k \alpha -\cos ^k \b |\leq \sqrt{k}\cdot \sin \gamma $$
 \end{itemize}
 \end{cor}

 \proof (i).
Recall   that if $M=(m_{ij})$ is
a $3\times 3$ positive semidefinite matrix, then   $(|m_{ij}|)$ is also positive semidefinite.
The first inequality in  (i)  is immediate from (\ref{one}). From a result of FitzGerald and Horn \cite[Theorem 2.2]{FitHor77}, we know  that
if $P=(p_{ij})$ is $3\times 3 $ nonnegative positive semidefinite, then the Hadamard power matrix
 $(p_{ij}^{r})$ is positive semidefinite for  all real $r\in [1, \infty )$. Thus
 for the $3\times 3$ positive semidefinite matrix  $B=(b_{ij})$ and for any real $k\geq 2$, matrix
 $(|b_{ij}|^{k/2})$ is  positive semidefinite. An application of (\ref{Eq:13xyz}) implies
  $$\big | |a|^k-|b|^k \big |=
  \big | (|a|^{k/2})^2-(|b|^{k/2})^2 \big |\leq \sqrt{1-(|c|^{k/2})^2}= \sqrt{1-|c|^{k}} $$

 (ii).  Let
 $a=\cos \a, b=\cos \b, c=\cos \g$. By
 Corollary \ref{cor3.14} (ii),  the matrix $B$ is positive semidefinite; so is the Hadamard power matrix $(b_{ij}^k)$. For $k=1$, $|a-b|\leq \sqrt{1-c^2}$ is  the same as $|\cos \a -\cos \b |\leq \sin \g$. For integer $k>1$,
 $$|a^k-b^k|\leq \sqrt{1-c^{2k}}=\sqrt{1-c^{2}} \cdot \sqrt{1+c^{2}+\cdots +c^{2(k-1)}}  \leq \sqrt{k}\cdot \sqrt{1-c^2}$$
which is the same as the desired inequality. $\qed$

In the following correlation matrix $D$,  $\sqrt{1-c}=0.3<|a-b|=0.4$. This
shows that the $\sqrt{2}$ in Corollary \ref{cor1ref}~(i) cannot be replaced by 1
in general.
$$D=\left ( \begin{array}{rrr}
  1  & 0    & 0.4  \\
  0  & 1    & 0.91 \\
0.4  & 0.91  & 1
 \end{array} \right )$$

In the previous proof, we saw $|\cos \a -\cos \b |\leq \sin \g$.
Inequality (\ref{Eq:13xyz}) implies $|\cos^2 \a -\cos^2 \b|\leq \sin \g$.
A question arises:  Can the $\sqrt{k}$ in Corollary \ref{cor1ref}~(ii) be removed or
replaced by a constant (like $\!\sqrt{2}\,$) that is independent of $k$?
The following result gives a negative answer.

\begin{thm}
Let  $0\leq  \a, \b, \g \leq \frac{\pi}{2}$, $|\a-\b |\leq \g \leq \a+\b$ and let
$$R_k=R_k(\a, \b, \g)=\left |\frac{\cos^k\a -\cos^k \b}{\sin \g} \right |, \quad \g \neq 0$$
Then for sufficiently large positive integer  $k$,
$$\sup_{\a, \b, \g}R_k \approx \sqrt{\frac{k}{e}}$$
\end{thm}

\proof If $\alpha =\beta$ then $R_k=0$. In addition, $R_k$ is symmetric with respect to $\alpha$ and $\beta$.
So, for the maximum $R_k$, we may assume $\alpha >\beta$. Moreover,
$$\frac{\cos^k\b -\cos^k \a}{\sin \g} \leq \frac{\cos^k\b -\cos^k \a}{\sin (\alpha-\beta)}$$

Let $R_k(\a, \b)=\frac{\cos^k\b -\cos^k \a}{\sin (\alpha-\beta)}$ if $\alpha \neq \b$ and
$R_k(\a, \b)= k\sin \b \cos^{k-1}\b$ if $\a =\b$ (denoted by $R_k(\b)$ for short).  Consider $R_k(\a, \b)$ over the
triangular region ${\Bbb D}=\{(\a, \b)\mid \frac{\pi}{2}\geq \a\geq \b \geq 0\}$.
 For every given  integer $k\geq 2$, the function
$R_k(\a, \b)$ is continuous on ${\Bbb D}$. To see this, it is sufficient to notice that
\begin{eqnarray*}
\lim_{(\a, \b)\rightarrow (t,t)}R_k(\a, \b) & = & \lim_{(\a, \b)\rightarrow (t,t)} \frac{(\cos \b -\cos \a)\big ( \sum_{l=0}^{k-1}\cos ^{k-1-l}\b \cos ^l\a\big )}{\sin (\a -\b)} \\
& = & \lim_{(\a, \b)\rightarrow (t,t)} \frac{\sin (\frac{\b + \a}{2})\big ( \sum_{l=0}^{k-1}\cos ^{k-1-l}\b \cos ^l\a\big )}{\cos (\frac{\a - \b}{2})} \\
 & = & k\sin t \cos^{k-1}t=R_{k}(t), \quad \mbox{where $t\in [0, \frac{\pi}{2}]$}.
\end{eqnarray*}
Thus $R_k(\a, \b)$ attains its maximum value at some point(s) in  ${\Bbb D}$.
On the $\a=\b$ portion of the boundary of $\Bbb D$, by computing the critical number,  the function $R_k(t)=k\sin t \cos^{k-1} t$ is maximized when $\tan t =\frac{1}{\sqrt{k-1}}$, and we get
$$\sup_{t\in [0, \frac{\pi}{2}]}R_k(t)={\frac{k}{\sqrt{k-1}}} \left (1-\frac1k \right )^{k/2} \;\; (\mbox{which is unbounded as $k\rightarrow \infty$})$$

For any given (small) $\varepsilon >0$, if $\a -\b\geq \varepsilon$, then
 $R_k(\a, \b)\leq \frac{1}{\sin \varepsilon}$ uniformly for all $k$.
So, in view of the behavior of $R_k(\a, \b)$ for $\a=\b$ in the above discussion, for all $k$ large enough,
  $R_k(\a, \b)$ is maximized as $\a$ approaches $\b$.
That is, if we let   $\a=\b +\varepsilon$, $\varepsilon>0$, then
\begin{eqnarray*}
\sup_{\a, \b, \g}R_k(\a, \b, \g) & = & \sup_{\a, \b}R_k (\a, \b)
        =\sup_{\b}\lim_{\varepsilon\rightarrow 0} R_k(\b+ \varepsilon, \b)\\
       & = & \sup_{\b} R_k(\b)=
        {\frac{k}{\sqrt{k-1}}} \left (1-\frac1k \right )^{k/2} \\
        & \approx & \sqrt{\frac{k}{e}} \quad \mbox{for sufficiently large $k$} \quad \qed
        \end{eqnarray*}

{\em Remark.}  In fact, one may show through some routine calculus computation that   for each $R_k(\a, \b)$ with $k\geq 1$, there are no interior critical points in $\Bbb D$.

\section{Metric-preserving and subadditive functions}

Let $a, b, c$ be nonnegative numbers.
The triplet $(a, b, c)\in \Bbb R^3$ is said to be  {\em triangle} if
$a\leq b+c,\; b\leq a+c,$ and $ c\leq a+b$; equivalently, $|a-b|\leq c\leq a+b$.
This can be restated as
$x\leq y+z$, where $x, y, z$ is any arrangement of $a, b, c$; equivalently,
   $|x-y|\leq z$ for any arrangement of $a, b, c$.
From the inequalities (t) and (T) of Section 2, we see for any nonzero vectors $u, v, w$ in an inner product space,
$\big (\theta (u, v), \theta (v, w), \theta (w, u)\big )$ and $\big (\Theta (u, v), \Theta (v, w), \Theta (w, u)\big )$
are   triangle triplets. They are special members of  the set of all triangle triplets:
$$\Delta  =\{ (a, b, c) \mid a, b, c\geq 0, |a-b|\leq c\leq a+b\}$$

For $a, b, c$  in $[0, 1]$, by Corollary \ref{cor3.14} (ii), the triplet $(\arccos a, \arccos b, \arccos c)$
is triangle if and only if
the $3\times 3$ matrix $B$ in Proposition \ref{Thm1:arccos}
 is positive semidefinite. Proposition \ref{Thm1:arccos} reveals a relation between the positivity of
 the $3\times 3$ matrix $B$ and triangle triplets
via certain functions.
This section is to present a theorem of this type, with which we show some inequalities for unit vectors.

Let $f$ be a nonnegative function defined on $[0, \infty )$. We say that $f$ is
{\em metric-preserving} (a metric preserver) if $f\circ d$ is also a metric on $M$, where $(M, d)$ is any metric space,
{\em triangle-preserving} (a triangle preserver) if
$(f(a), f(b), f(c))$ is triangle whenever $(a, b, c)$ is triangle, and
{\em subadditive} if $f(s+t)\leq f(s)+f(t)$
for all $s, t\geq 0$.
The reader is referred to \cite{Cor99} 
for metric-preserving functions
and  \cite[Chapter 16]{Kuc09} and \cite[Chapter 12]{Sch96})  for subadditive functions.

These three functions are closely related, but not exactly the same.
It is known (see, e.g., \cite[p.\,9]{Dob1998}) that
if $f:  [0, \infty )\mapsto [0, \infty )$ is nondecreasing
and subadditive then $f$ is a triangle preserver.
A nonnegative concave function vanishing at 0 is necessarily subadditive (see, e.g., \cite[p.\,314]{Sch96}).
Nonnegative concave (not necessarily continuous) functions must be nondecreasing. (This seems to be a known fact; but we
were not able to find a reference with a proof.)  A stronger
version of the result is stated as:

\begin{pro}
Let $L=[l, \infty )$ or $(l, \infty)\subseteq \Bbb R$. If $f$ is   nonnegative   on
$L$, i.e.,  $f(x)\geq 0$ for all $x\in L$,    and if $f$ is mid-point concave on $L$, i.e.,
$$ f\Big (\frac12 x+\frac12 y\Big )\geq \frac12 f(x)+\frac12 f(y),
\;\;\mbox{for all $x, y\in L$}$$
then $f(x)$ is monotonically increasing on $L$, i.e., $f(x)\geq f(y)$ for $x>y$, $x, y\in L$.
\end{pro}

\proof
Suppose that $f(x)$ is not monotonically increasing. Then there exist $s, t\in L$, $t>0$,  such that $f(s)>f(s+t)$.
Let $r=f(s)-f(s+t)>0$.
Since
$$f(s+t)=f\Big ( \frac12 s+ \frac12 (s+2t)\Big )\geq \frac12 f(s)+\frac12 f(s+2t)$$
we arrive at
$$f(s+t)-f(s+2t)\geq f(s)-f(s+t)=r$$
Let
$$F_n=f(s+nt)-f\big (s+(n+1)t\big ), \quad n=0, 1, 2, \dots$$
Then, in a similar way as above, we can show that $\{F_n\}$ is a decreasing sequence bounded by $r$ from below,
i.e., $F_{n}\geq F_{n-1}\geq \cdots \geq F_1\geq r$. It follows that
$$f(s)-f(s+nt)=
F_0+F_1+F_2+\cdots + F_{n-1}\geq nr.$$
Thus, $f(s+nt)<0$   when $n$ is large enough,  contradicting $f(x)\geq 0$, $x\in L$.  $\qed$

\begin{thm}\label{Thm:140425}
Let $a, b, c$ be real numbers such that the $3\times 3$ matrix
$$B=\left ( \begin{array}{ccc}
  1 & a & b  \\
  a  & 1 & c \\
b & c   & 1
 \end{array} \right )$$
 is positive semidefinite. Then for all functions $f$ described
 in Proposition \ref{Thm1:arccos} and
 for all  metric-preserving functions $g$, with  $h=g\circ f^{-1}$, we have
 $$h(a)\leq  h(b)+h(c)$$
 Consequently, for any real $k\geq 2$,
 $$\sqrt[k]{1-|a|^k}\leq \sqrt[k]{1-|b|^k}+\sqrt[k]{1-|c|^k}$$
 \end{thm}

\proof Proposition \ref{Thm1:arccos} says that $(f^{-1}(a), f^{-1}(b), f^{-1}(c))$ is a triangle triplet. For any
metric-preserving function $g$, with $h=g\circ f^{-1}$,
$(h(a), h(b), h(c))$ is also a triangle triplet.
This gives the desired inequality.

 For the second part, for any fixed real $k\geq 2$, consider the function
 $$p(t)=\left \{ \begin{array}{cl}
 \sqrt[k]{1-\cos ^k t} &  \mbox{if $0\leq t\leq \frac{\pi}{2}$} \\
 1 & \mbox{if $t>\frac{\pi}{2}$} \end{array} \right .
 $$

 It is straightforward to verify that $p(t)$ is nonnegative, increasing,  and
 concave for $k\geq 2$ (by checking $p'(t)\geq 0$, $p''(t)\leq 0$); thus, $p(t)$ is
 metric-preserving.
 Applying $p(t)$ to $(\arccos |a|,  \arccos |b|, \arccos |c|)$ yields the inequality. $\qed$

 \begin{cor} Let $u, v, w$ be any unit vectors of an inner product space. Then
 \begin{equation}\label{Ineq:lin}
 \sqrt[k]{1-|\langle u, v\rangle |^k}\leq \sqrt[k]{1-|\langle u, w\rangle |^k}+\sqrt[k]{1-|\langle w, v\rangle |^k}
 \end{equation}
 for any real number $k\geq 2$. In particular,
 \begin{equation}\label{sqrt||}
 \sqrt{1-|\langle u, v\rangle |^2}\leq \sqrt{1-|\langle u, w\rangle |^2}+\sqrt{1-|\langle w, v\rangle |^2}
 \end{equation}
Similar inequalities hold for $|\re \langle \cdot  , \cdot \rangle |$ in place of $|\langle \cdot  , \cdot \rangle |$.
\end{cor}

Inequality (\ref{sqrt||}) appears in \cite{WZMonthly} (see also \cite[p.\,195]{ZFZbook11}). Inequality (\ref{Ineq:lin})
is seen in \cite{Lin12_MathIntel} (with a minor glitch on the condition $b+c\leq 1$ which can be fixed).

\begin{cor}\label{dc}
Let $u, v, w$ be unit vectors in an inner product space,  $\a, \b, \g$ be respectively
the angles  $\theta (u, v), \theta (v, w), \theta (w, u)$ or $\Theta (u, v), \Theta (v, w), \Theta (w, u)$.
Then
\begin{itemize}
\item[\rm (i).] $\a \leq \b +\g$.
\item[\rm (ii).] $\sin \a \leq \sin \b + \sin \g$.
\item[\rm (iii).]  $\cos \a \not \leq \cos \b + \cos \g$   in general.
\item[\rm (iv).]  $\cos \a \leq \cos \b + \sin \g$.
\end{itemize}
\end{cor}

\proof
(i) is the same as (t) and (T) in Section 2. (ii) is true because of (\ref{sqrt||}) and
the similar inequality for $\re \langle \cdot  , \cdot \rangle $.
For (iii), take
 $u=(0, 0, 1)$, $v=\frac{1}{\sqrt{2}}
(1, 0, 1)$,  and $w=(0, 1, 0)$ in $\mathbb{R}^{3}$   with the standard Euclidean inner
product. (iv) follows from inequality (\ref{one}).
$\qed$

 \medskip

\noindent
{\bf Acknowledgement} The authors are thankful to the referee for several comments and
suggestions that improved the results and exposition of the paper.

\end{document}